\documentclass[12pt]{amsart}
\usepackage{amssymb, amstext, amscd, amsmath}
\makeatletter
\def\@cite#1#2{{\m@th\upshape\bfseries%
[{#1\if@tempswa{\m@th\upshape\mdseries, #2}\fi}]}} \makeatother
\theoremstyle{plain}
\newtheorem{thm}{Theorem}[section]
\newtheorem{cor}[thm]{Corollary}
\newtheorem{prop}[thm]{Proposition}
\newtheorem{lem}[thm]{Lemma}
\theoremstyle{definition}
\newtheorem{rem}[thm]{Remark}
\newtheorem{defn}[thm]{Definition}

\newtheorem{eg}[thm]{Example}
\newcommand{\Prf}{\noindent\textbf{Proof.\ }}
\newcommand{\bx}{\hfill$\blacksquare$\medbreak}

\newcommand{\ca}{\mathrm{C}^*}

\newcommand{\ol}{\overline}
\newcommand{\td}{\widetilde}

\newcommand{\sot}{\textsc{sot}}
\newcommand{\wot}{\textsc{wot}}

\newcommand{\bbB}{{\mathbb{B}}}
\newcommand{\bbC}{{\mathbb{C}}}
\newcommand{\bbF}{{\mathbb{F}}}

\newcommand{\bbN}{{\mathbb{N}}}

\newcommand{\bbT}{{\mathbb{T}}}
\newcommand{\bbZ}{{\mathbb{Z}}}

  \newcommand{\A}{{\mathcal{A}}}
  \newcommand{\B}{{\mathcal{B}}}

  \newcommand{\E}{{\mathcal{E}}}
  \newcommand{\F}{{\mathcal{F}}}
  \newcommand{\G}{{\mathcal{G}}}
\renewcommand{\H}{{\mathcal{H}}}
  
  \newcommand{\J}{{\mathcal{J}}}

  \newcommand{\M}{{\mathcal{M}}}
  \newcommand{\N}{{\mathcal{N}}}

\newcommand{\fA}{{\mathfrak{A}}}

\newcommand{\fL}{{\mathfrak{L}}}

\newcommand{\fR}{{\mathfrak{R}}}

\newcommand{\qand}{\quad\text{and}\quad}

\newcommand{\qfor}{\quad\text{for}\quad}

\newcommand{\Alg}{\operatorname{Alg}}

\newcommand{\Aut}{\operatorname{Aut}}

\newcommand{\dist}{\operatorname{dist}}

\newcommand{\Lat}{\operatorname{Lat}}

\newcommand{\rad}{\operatorname{rad}}
\newcommand{\nc}{\operatorname{NC}}

\newcommand{\spn}{\operatorname{span}}

\newcommand{\obj}{\operatorname{Obj}}

\newcommand{\flgee}{\fL_G}

\newcommand{\flamd}{\fL_{(\Lambda,d)}}

\begin{document}

%%%%%%%%%%%%%%%%%%%%%%%%%%%%%%%%%%%%%%%%%%%%%%
%%%%%%%%%%%%%%
\title[$H^\infty$ Algebras of Higher Rank Graphs]{The $H^\infty$ Algebras of Higher Rank Graphs}
%\thanks{}
%
%
\author[D.W. Kribs]{David~W.~Kribs}
\address{Department of Mathematics and Statistics, University of
Guelph, Guelph, Ontario Canada N1G 2W1} \email{dkribs@uoguelph.ca}

\author[S.C. Power]{Stephen~C.~Power}
\address{Department of Mathematics and Statistics, Lancaster
University, Lancaster, United Kingdom  LA1 4YF}
\email{s.power@lancaster.ac.uk}
\begin{abstract}
We begin the study of a new class of operator algebras that arise
from higher rank graphs. Every higher rank graph generates a Fock
space Hilbert space and creation operators which are partial
isometries acting on the space. We call the weak operator topology
closed algebra generated by these operators a {\it higher rank
semigroupoid algebra}. A number of examples are discussed in
detail, including the single vertex case and higher rank cycle
graphs. In particular the cycle graph algebras are identified as
matricial multivariable function algebras. We obtain reflexivity
for a wide class of graphs and characterize semisimplicity in
terms of the underlying graph.
\end{abstract}

\date{}
\maketitle
%%%%%%%%%%%%%%%%%%%%%%%%%%%%%%%%%%%%%%%%%%%%%%
%%%%%%%%%%%%%%

%%%%%%%%%%%%%%%%%%%%%%%%%%%%%%%%%%%%%%%%%%%%%%%%%%%
%\section{Introduction}\label{S:intro}
%%%%%%%%%%%%%%%%%%%%%%%%%%%%%%%%%%%%%%%%%%%%%%%%%%%

In \cite{KumP3} Kumjian and Pask  introduced $k$-graphs as an
abstraction of the combinatorial structure underlying the higher
rank graph $\ca$-algebras of Robertson and Steger \cite{RS2,RS1}.
A $k$-graph generalizes the set of finite paths of a countable
directed graph when viewed as a partly defined multiplicative
semigroup with vertices considered as degenerate paths. The
$\ca$-algebras associated with $k$-graphs include $k$-fold tensor
products of graph $\ca$-algebras, and much more
\cite{APS,KumP4,PQR2,PQR1,RaSi}. On the other hand, as a
generalization of the nonselfadjoint free semigroup algebras
$\fL_n$ \cite{ArPop,DKP,DP1,DP2,Kfactor,Pop1,Pop4}, the authors
\cite{KP1,KP2} have recently studied free semigroupoid algebras
$\fL_G$ associated with  directed countable graphs $G$. In
particular it was shown that these algebras are reflexive. (See
also \cite{JP,JK2,JK1,KK1,KK2,Mu,MS1,MS2,Sol} for related recent
work.) As it turns out, these algebras arise from the left regular
representation of the 1-graph of the directed graph $G$. In the
present paper we consider the higher rank versions of these
algebras, the $k$-graph algebras $\flamd$ associated with the
$k$-graph $(\Lambda,d)$, as well as their norm closed subalgebras.
To our knowledge such nonselfadjoint higher rank graph algebras
have not been considered previously. However, from the perspective
of contemporary operator algebra theory they evidently form a
natural class and one which may play an important role in more
general higher rank operator algebra considerations.

The algebras $\fL_n$ are also referred to as the noncommutative
analytic Toeplitz algebras and in the case $n=1$ one obtains the
usual algebra $H^\infty$ acting on the Hardy space of the circle
\cite{Douglas_text,Hoffman_text,Sarason}. In the present paper we
examine eigenvalues, reflexivity, hyper-reflexivity and
semisimplicity for the algebras $\fL_{(\Lambda,d)}$ which can be
viewed as the $H^\infty$ algebras of higher rank graphs.

In  $\S 1$ we outline the nomenclature associated with higher rank
graphs $(\Lambda,d)$. In $\S 2$ we introduce higher rank
semigroupoid algebras $\fL_{(\Lambda,d)}$ and derive some basic
properties. We follow this in $\S 3$ by presenting a diverse
collection of examples and in $\S 4$ we consider the single vertex
algebras. In particular we determine the eigenvalues for the
adjoint algebras and the Gelfand space of the norm closed
subalgebras $\A_{(\Lambda,d)}$. In the next section ($\S 5$) we
prove the algebras $\fL_{(\Lambda,d)}$ are reflexive or
hyper-reflexive for various diverse graphs. In the final section
($\S 6$) we find a graph condition which characterizes when
$\fL_{(\Lambda,d)}$ is semisimple and give an explicit description
of the Jacobson radical in the finite vertex case.

%%%%%%%%%%%%%%%%%%%%%%%%%%%%%%%%%%%%%%%%%%%%%%%%%%%
\section{Higher Rank Graphs}\label{S:kgraphs}
%%%%%%%%%%%%%%%%%%%%%%%%%%%%%%%%%%%%%%%%%%%%%%%%%%%

Let $G = (V,E)$ be a countable directed graph and let $\Lambda_G$
denote the set of all directed paths $\lambda = e_r e_{r-1} \cdots
e_1$ where $e_1$ is an edge $(v_2,v_1)$ directed from $v_1$ to
$v_2$ and where $e_k$ is an edge $(v_{k+1},v_{k})$, for $1\leq k
\leq r$, where $v_1, v_2, \ldots, v_r$ are the vertices of the
path, in their directed order, possibly with repetitions. Let $d:
\Lambda_G \rightarrow \bbN$ be the length function. Then the pair
$(\Lambda_G, d)$ is an example of a 1-graph in the sense of the
definition below.

The set $\Lambda_G$ has a natural partially defined multiplication
which together with vertices as degenerate paths makes $\Lambda_G
\cup V$ into a  (discrete) {\it semigroupoid} with vertices as
units. (This is the terminology of \cite{KP1,KP2}.) However,
$\Lambda_G$ can be viewed as a set of morphisms between elements
of $V$, and as such $\Lambda_G$ forms a small category with $V$ as
the set of objects. (`Small' since the objects form a set.) It is
this viewpoint which is extended in the definition below.

\begin{defn}\label{hrdefn}
\cite{KumP3} A $k$-graph $(\Lambda, d)$ consists of a countable
small category $\Lambda$, with range and source maps $r$ and $s$
respectively, together with a functor $d: \Lambda \rightarrow
\bbZ_+^k$ satisfying the factorization property: for every
$\lambda\in\Lambda$ and $m,n\in\bbZ_+^k$ with $d(\lambda) = m+n$,
there are unique elements $\mu,\nu\in\Lambda$ such that $\lambda =
\mu\nu$ and $d(\mu) = m$ and $d(\nu) =n$.
\end{defn}

By the factorization property we may identify the objects
$\obj(\Lambda)$ of $\Lambda$ with the subset $\Lambda^0 \equiv
d^{-1} (0,\ldots ,0)$. We also write $\Lambda^n$ for $d^{-1}(n)$,
$n\in \bbZ_+^k$. Observe that the factorization property implies
that left and right cancellation hold in $\Lambda$. Let $\delta :
\Lambda \rightarrow \bbZ_+$ be the grading function defined by
$
\delta(\lambda) = |d(\lambda)| = n_1 + \ldots + n_k
$
where $d(\lambda) = (n_1,\ldots ,n_k)$.

As indicated above the conventional definition of a directed
graph, with its semigroupoid of paths, is captured in the special
case of 1-graphs. The set $\Lambda^n$ corresponds to the directed
paths of length $n$ in the graph. We give a variety of examples of
$k$-graphs in $\S$~\ref{S:examples}. However, let us consider here
the simple example given by the 2-graph $\Lambda_{(G_1 \times G_2,
d)}$ arising from the set $\Lambda_{G_1\times G_2}$ of paths in
the direct product directed graph $G_1\times G_2$, where $G_1$,
$G_2$ are directed graphs, together with the natural map $d:
\Lambda_{G_1\times G_2} \rightarrow \bbZ_+^2$. In this case, if $e
= (v_2, v_1)$ is an edge of $G_1$ and $f=(w_2,w_1)$ is an edge of
$G_2$, then $e\times f = ((v_2,w_2),(v_1,w_1))$ is an edge of
$G_1\times G_2$ and $d(e\times f) = (1,1)$. Also, by definition
$G_1\times G_2$ includes all edges $((v,w_2),(v,w_1))$ and
$((v_2,w),(v_1,w))$, with $d$-degrees $(0,1)$ and $(1,0)$
respectively, for each vertex $v$ of $G_1$ and $w$ of $G_2$. More
generally, a direct product $G_1\times \ldots \times G_k$ of $k$
directed graphs generates a $k$-graph in this way.

It is of immediate interest to see $k$-graphs which do not arise
as direct products. For an elementary example let $\Lambda^0 = \{
v\}$, $\Lambda^{(1,0)} = \{a\}$, $\Lambda^{(0,1)} = \{ b\}$ be
singleton sets. Suppose moreover that composition of the morphisms
$a,b$ generate all morphisms of the category $\Lambda$. In view of
the required factorization property and its uniqueness, it soon
becomes clear that all the morphisms $w= a^{n_1} b^{m_1} a^{n_2}
\cdots b^{m_r}$ with degree defined by $d(w) = (n,m)$,
$n=n_1+\ldots + n_r$, $m= m_1+ \ldots + m_r$ must coincide. Thus
set $\Lambda^{(n,m)} =\{a^n b^m\}$ and in this way we obtain a
2-graph $\Lambda = \cup_{n\in\bbZ_+^2} \Lambda^n$.

Note that the 2-graph above is generated by the units and elements
of total degree 1 subject to a simple commutation relation. We now
give a similar such description of more general 2-graphs in which
a commutation rule $\alpha \times \beta =
\theta(\alpha\times\beta)$ is built in to ensure the factorization
property.

Let $A = \Lambda_{(G_1,d_1)}$ and $B = \Lambda_{(G_2,d_2)}$ be
1-graphs such that $A^0 = B^0$, so that the underlying graphs
$G_1$ and $G_2$ have the same number of vertices and the vertex
sets are identified. Let $v_1$, $v_2$ be two vertices and consider
the following sets of pairs of edges,
\begin{align*}
E(v_2,v_1) &= \Big\{ \alpha\times\beta\in A^1 \times B^1 \,\,\Big|
\,\,s(\alpha)=r(\beta), \, s(\beta) = v_1, \, r(\alpha) =v_2\Big\} \\
F(v_2,v_1) &= \Big\{ \beta\times\alpha\in B^1 \times A^1
\,\,\Big|\,\, s(\beta)=r(\alpha), \, s(\alpha) = v_1, \, r(\beta)
=v_2\Big\}.
\end{align*}
Suppose also that these sets have the same cardinality for all
vertex pairs and that $\theta$ is a bijection mapping each
$\alpha\times\beta$ in $E(v_2,v_1)$ to an element
$\theta(\alpha\times\beta)$ in $F(v_2,v_1)$, for all vertex pairs.
To construct the 2-graph $(\Lambda,d) = A \ast_\theta B$ define
\[
\Lambda^0 = V = V(G_1) = V(G_2),
\]
\[
\Lambda^{(1,0)} = \big\{ \alpha \times v \in A^1 \times V
\,\,\big| \,\, s(\alpha) = v\big\},
\]
\[
\Lambda^{(0,1)} = \big\{ v \times \beta \in V \times B^1 \,\,\big|
\,\, r(\beta) = v\big\},
\]
\[
\Lambda^{(1,1)} = \bigcup_{v,w\in V} E(v,w) = \bigcup_{v,w\in V}
F(v,w),
\]
where the last equality arises from the identifications
$\alpha\times\beta = \theta(\alpha\times\beta)$. Plainly the
factorization property holds for morphisms in $\Lambda^{(1,1)}$.
Finally define $\Lambda^{(n,m)}$ as the set of morphisms obtained
from arbitrary finite compositions of morphisms in
$\Lambda^{(1,0)}$ and $\Lambda^{(0,1)}$ subject to the relations
generated by the identifications $(\alpha\times v)(v\times\beta) =
(\beta_1 \times v)(v\times \alpha_1)$ if $\alpha_1\times \beta_1 =
\theta(\alpha\times\beta)$. It is routine to check that
$\Lambda^{(n,m)}$ satisfies the factorization property  with the
natural map $d:\Lambda \rightarrow \bbZ_+^2$ (where $\Lambda$ is
the union of all the sets $\Lambda^{(n,m)}$), and thus the pair
$(\Lambda, d)$ is a 2-graph.

\begin{rem}
We shall find it convenient to view a 2-graph as being specified
through a directed graph in which edges are of two types, perhaps
red or blue, according to their degree, $(1,0)$ or $(0,1)$,
together with a set of relations that define the factorization
property.  For instance, if $\lambda = e_3e_2e_1$ is a path with
$d(\lambda)=(2,1)$, $\delta(\lambda)=3$ and $e_1,e_2,e_3$ coloured
red, red, blue respectively, then there must be red edges
$f_1,f_3,g_2,g_3$ and blue edges $f_2,g_1$ such that $\lambda =
f_3f_2f_1 = g_3g_2g_1$. The initial (respectively final) vertices
of $e_1,f_1,g_1$ (respectively $e_3,f_3,g_3$) must coincide  but
$e_2,f_2,g_2$ may have distinct initial and final vertices.

It is thus understood that a given path in the chromatic graph
represents an equivalence class of paths under the commutation
relations which are either specified explicitly, or implied by the
factorization property. Of course, the same remarks apply to a
$k$-graph, which corresponds to a $k$-coloured graph
$\Lambda^{(e_1)} \cup \ldots \cup \Lambda^{(e_k)}$ together with
commutation relations distinct colours.
\end{rem}

\begin{rem}
Let us clarify our use of the terminology ``semigroupoid'' and
``freeness''. To each directed graph $G$ one can associate the
(universal) graph $\ca$-algebra, and under mild hypotheses this is
isomorphic to a (topological) groupoid $\ca$-algebra $\ca(\G)$
associated with the topological path groupoid $\G$ of $G$. We have
no cause in this paper to consider this groupoid but we do find it
convenient to use terminology which derives from the (discrete)
groupoid of an undirected graph $G$. This consists of paths in the
edges $e$ of $G$ and their formal inverses $e^{-1}$, together with
the vertices viewed as degenerate edges forming units. With the
understanding that the only identification of paths is through the
relations $ee^{-1} = r(e)$ and $e^{-1}e = s(e)$, it is natural to
refer to this groupoid as the {\it free} groupoid $\bbF(G)$ of
$G$. Indeed, in the case of a single vertex graph this groupoid is
the free group on $n$-generators where $n$ is the number of edges
of $G$. Moreover, if $\E$ is a discrete groupoid generated by
elements $e_1, e_2, \ldots$ together with units $x_1, x_2,
\ldots$, then this set of generators determines a graph, $G$ say.
If $\alpha : \E \rightarrow \F$ is a discrete groupoid
homomorphism then there is a lifting $\beta : \bbF(G) \rightarrow
\F$ such that $\beta = \pi \circ\alpha$ where $\pi : \bbF(G)
\rightarrow \E$ is the natural map. Thus, $\bbF(G)$ is the free
object in the category of discrete groupoids with generators
labelled by the graph $G$.

Similarly, if we omit the formal inverses $e^{-1}$ of the edges of
a graph $G$ then we identify a unital semigroupoid in $\bbF(G)$
which we denote as $\bbF^+(G)$ and refer to as the (discrete) free
semigroupoid of $G$. Thus, a 1-graph coincides with the pair
$(\bbF^+(G),d)$ where $d$ is the length function and $G$ is the
graph arising from elements of total degree 1.
\end{rem}

%%%%%%%%%%%%%%%%%%%%%%%%%%%%%%%%%%%%%%%%%%%%%%%%%%%%%%%%%%%%%%%
\section{Higher Rank Semigroupoid Algebras}\label{S:hrsa}
%%%%%%%%%%%%%%%%%%%%%%%%%%%%%%%%%%%%%%%%%%%%%%%%%%%%%%%%%%%%%%%

Let $(\Lambda,d)$ be a $k$-graph. Let $\H_\Lambda$ be the Fock
space of $\Lambda$ which we define to be the Hilbert space with
orthonormal basis $\{\xi_\lambda : \lambda \in \Lambda\}$. For
$\lambda\in\Lambda$ define the operator $L_\lambda$ on
$\H_\lambda$ such that
\[
L_\lambda \xi_\mu = \left\{
\begin{array}{cl}
\xi_{\lambda\mu} & \mbox{if $s(\lambda)=r(\mu)$} \\
0 & \mbox{if $s(\lambda)\neq r(\mu)$}
\end{array} \right..
\]
It follows from the factorization property that each $L_\lambda$
is a partial isometry. Moreover, $L_v$, $v\in\Lambda^0$, is the
projection onto the subspace $\spn\{\xi_\lambda: r(\lambda)=v\}$.

\begin{defn}
The {\it semigroupoid algebra} $\fL_{(\Lambda,d)}$ of the
$k$-graph $(\Lambda,d)$ is the weak operator topology closed
linear span of $\{L_\lambda : \lambda\in\Lambda\}$.
\end{defn}

Recalling the definition of the grading function $\delta$, it is
evident that $\H_\lambda$ is naturally graded as $\H_\lambda =
\H_0 \oplus \H_1 \oplus \H_2 \oplus \ldots$, where $\H_n$ is the
closed span of the $\xi_\lambda$ with $\delta(\lambda) = n$.

By arguing exactly as in the case of free semigroupoid algebras
\cite{KP1} one can obtain the following proposition. For brevity
we write $\fL_\Lambda$ for $\flamd$ and $\fR_\Lambda$ for the
analogue of $\fL_\Lambda$ for right actions.

\begin{prop}
If $A\in\fL_\Lambda$ then $A$ is the $\sot$-limit of the Cesaro
sums
\[
\sum_{\delta(\lambda)\leq n} \Big( 1 -\frac{\delta(\lambda)}{n}
\Big) a_\lambda L_\lambda,
\]
where $a_\lambda\in\bbC$ is the coefficient of $\xi_\lambda$ in
$A\xi_v = \sum_{s(\lambda) =v} a_\lambda \xi_\lambda$, for
$v\in\Lambda^0$.
\end{prop}

Thus elements of $\fL_\Lambda$ have Fourier expansions  $A \sim
\sum_{\lambda\in\Lambda} a_\lambda L_\lambda$. In particular, this
leads to the following description of the commutant.

\begin{prop}\label{commprop1}
The commutant of $\fR_\Lambda$ is $\fL_\Lambda$.
\end{prop}

\Prf As in the directed graph case we can consider the Cesaro
operators associated with the partition $I= E_0+E_1+\ldots$ where
$E_n$ is the projection onto $\H_n$. These operators are given by
\[
\Sigma_n(A) = \sum_{\delta(\lambda)=m<n} \Big(
1-\frac{\delta(\lambda)}{n}\Big) \Phi_m(A),
\]
where the operators $\Phi_m(A) = \sum_{n\geq \max\{ 0,-m\}} E_n A
E_{n+m}$ are the diagonals of $A$ with respect to the block
matrix
decomposition associated with the partition.
The operators $\Sigma_n(A)$ converges to $A$ in the strong
operator topology for all $A \in\B(\H_\Lambda)$.

It is clear that $\fL_\Lambda$ is contained in
$\fR_\Lambda^\prime$, thus for  the converse we fix $A \in
\fL_\Lambda^\prime$. We will show that $A_v \equiv AL_v$ belongs
to $\fL_\Lambda$ for all $v\in \Lambda^0$. This will finish the
proof since $A = \sum_{v\in\Lambda^0} AL_v$, the sum converging
$\sot$ when $\Lambda^0$ is infinite. Let $A\xi_v =R_v A_v\xi_v=
\sum_{s(\lambda)= v} a_\lambda \xi_\lambda.$ Define operators in
$\fL_\Lambda$ by
\[
p_n(A_v) = \sum_{ \delta(\lambda) <n ; \, s(\lambda)= v} \Big( 1 -
\frac{\delta(\lambda)}{n}\Big) a_\lambda L_\lambda.
\]
We will prove that $A_v = \sot \!-\!\lim_{n\rightarrow\infty}
p_n(A_v)$ by showing that $p_n(A_v) =\Sigma_n(A_v)$. First note
that $\Phi_m(A_v)$ belongs to $\fR_\Lambda^\prime$ for all $m$
since $A_v$ belongs to $\fR_\Lambda^\prime$ and $E_{n+1} R_\lambda
= R_\lambda E_n$ for all $n$ and $\lambda\in\Lambda^1$, while
$\Phi_m(A_v)$ commutes with each projection $R_w$ since $R_w E_n =
E_n R_w$ is the projection onto $\spn\{ \xi_\lambda:
\delta(\lambda) =n, \,\, s(\lambda)= w \}$ for all $n$. It follows
that $\Sigma_n(A_v)$ belongs to $\fR_\Lambda^\prime$ for $n\geq
1$.

Now it is enough to show that $\Sigma_n(A_v) \xi_v = p_n(A_v)
\xi_v$. If this is the case, then for $\lambda\in\Lambda$ with
$r(\lambda)=v$ we have
\[
\Sigma_n(A_v) \xi_\lambda = R_\lambda \Sigma_n(A_v) \xi_v =
R_\lambda p_n(A_v) \xi_v = p_n(A_v) \xi_\lambda,
\]
whereas, if $r(\lambda)=w$ with $w \neq v$ then
\[
\Sigma_n(A_v)\xi_\lambda = \Sigma_n(A_v)L_v L_w \xi_\lambda =0=
p_n(A_v)L_w \xi_\lambda = p_n(A_v)\xi_\lambda,
\]
since $L_v$ commutes with each $E_n$.

Observe that $\Phi_0 (A_v) \xi_v = E_0(A_v)E_0\xi_v = a_v \xi_v$,
and $\Phi_m(A_v)\xi_v = 0$ for $m>0$. Further, for $m< 0$ we have
\[
\Phi_m(A_v) \xi_v = (E_{-m} A_v )\xi_v = E_{-m} \sum_{s(\lambda)=
v} a_\lambda \xi_\lambda = \sum_{s(\lambda) v;
\,\delta(\lambda)=-m} a_\lambda \xi_\lambda.
\]
Hence it follows that
\[
\Sigma_n(A_v) \xi_v =   \sum_{\delta(\lambda) <n;\, s(\lambda)= v}
\Big( 1 -\frac{\delta(\lambda)}{n}\Big) a_\lambda \xi_\lambda =
p_n(A_v) \xi_v,
\]
as required. Therefore each $A_v = AL_v$ belongs to $\fL_\Lambda$
and this completes the proof. \bx

Given $\Lambda$, let $\Lambda^t$ be a category with the same
functor $d$, with $\obj(\Lambda^t) = \obj(\Lambda)$ and morphisms
for each $\lambda\in\Lambda$ denoted by $\lambda^t$ with
$s(\lambda^t) = r(\lambda)$ and $r(\lambda^t) = s(\lambda)$. Then
a simple argument shows that $\fL_\Lambda$ and $\fR_{\Lambda^t}$
are unitarily equivalent via the unitary $U: \H_{\Lambda^t}
\rightarrow \H_\Lambda$ defined by $U\xi_{\lambda^t} =
\xi_\lambda$.

\begin{cor}\label{leftcommutant}
The commutant of $\fL_\Lambda$ is $\fR_\Lambda$.
\end{cor}

\Prf If $U$ is the unitary above then $\fR_{\Lambda^t}^\prime =
(U^*\fL_\Lambda U)^\prime = U^* \fL_\Lambda^\prime U$. Hence by
Proposition~\ref{commprop1} we have $\fR_\Lambda =
U\fL_{\Lambda^t}U^*= U \fR_{\Lambda^t}^\prime U^* =
\fL_\Lambda^\prime. $ \bx

\begin{cor}
$\fL_\Lambda$ is its own second commutant, $\fL_\Lambda
=\fL_\Lambda^{\prime\prime}$.
\end{cor}

%%%%%%%%%%%%%%%%%%%%%%%%%%%%%%%%%%%%%%%%%%%%%%%%%%%
\section{Examples}\label{S:examples}
%%%%%%%%%%%%%%%%%%%%%%%%%%%%%%%%%%%%%%%%%%%%%%%%%%%

We now describe a number of examples of higher rank semigroupoid
algebras starting with some elementary direct product $k$-graphs.

\begin{eg}\label{1steg}
Let $C_1$ be the directed graph with a single vertex $v$ and loop
edge $e$.  Then the Fock space $\H_{\Lambda_{C_1}}$ may be
identified with the Hardy space $H^2$ and under this
identification $\fL_{\Lambda_{C_1}}$ is unitarily equivalent to
the analytic Toeplitz algebra $H^\infty$
\cite{Douglas_text,Hoffman_text,Sarason}. Consider, as in
$\S$~\ref{S:kgraphs},  the natural direct product $\Lambda =
\Lambda_{C_1}\times \Lambda_{C_1}$ and let $\ol{v} = v\times v$,
$a= e\times v$ and $b= v\times e$. Then $\Lambda^0 = \{ \ol{v}
\}$, $\Lambda^{(1,0)} = \{ a \}$, $\Lambda^{(0,1)} = \{ b \}$ and
it becomes clear that $\Lambda$ is the simple 2-graph discussed in
$\S$~\ref{S:kgraphs}. The standard basis for the Fock space
$\H_\Lambda \cong H^2 \otimes H^2$ may be identified with the
vertices in the 2-lattice of positive integers $\bbZ_+^2$ and
$\fL_\Lambda$ is unitarily equivalent to $H^\infty \otimes
H^\infty$.

More generally, given directed graphs $G_1, \ldots , G_k$  the
standard basis for the Fock space $\H_{\Lambda_{G_1}\times \ldots\times \Lambda_{G_k}}$ may be identified with the standard basis
for $\H_{G_1}\otimes \ldots \otimes \H_{G_k}$ and this
identification yields the unitary equivalence
$\fL_{\Lambda_{G_1}\times \ldots \times \Lambda_{G_k}} \cong
\fL_{G_1} \otimes \ldots \otimes \fL_{G_k}$. For example, if $F_n$
is the directed graph with a single vertex and $n\geq 2$ distinct
loop edges, then $\fL_{F_n} = \fL_n$ is the free semigroup algebra
(the noncommutative analytic Toeplitz algebra) which acts on
unrestricted $n$-variable Fock space $\H_{F_n} \equiv \H_n$. The
standard basis for $\H_n$ is identified with the set of all words
from an alphabet with $n$ noncommuting letters. Thus, given
positive integers $n_1, \ldots ,n_k$ we have
$\H_{\Lambda_{F_{n_1}} \times \ldots \times \Lambda_{F_{n_k}}}
\cong \H_{n_1} \otimes \ldots \otimes \H_{n_k}$ and
$\fL_{\Lambda_{F_{n_1}} \times \ldots \times \Lambda_{F_{n_k}}}
\cong \fL_{n_1} \otimes \ldots \otimes \fL_{n_k}$.
\end{eg}

\begin{eg}
Let $G$ be the connected directed graph with two edges $a_1 =
(x_2,x_1)$, $a_2 = (x_3,x_2)$. Then the free semigroupoid algebra
$\flgee$ is unitarily equivalent to the operator algebra of
matrices
\[
\left[ \begin{matrix} \alpha & 0 & 0 \\ \delta & \beta & 0 \\
\kappa & \epsilon & \gamma
\end{matrix} \right] \quad \oplus \quad \left[ \begin{matrix} \beta & 0  \\ \epsilon &
\gamma
\end{matrix} \right] \quad \oplus \quad \left[ \begin{matrix}
\gamma
\end{matrix} \right],
\]
where $\alpha,\beta,\gamma,\delta,\epsilon,\kappa\in\bbC$, acting
on the Fock space
\[
\H_G = \big( \bbC \xi_{x_1} + \bbC \xi_{a_1} + \bbC \xi_{a_2a_1}
\big)\,\, \oplus \,\,\big( \bbC \xi_{x_2} + \bbC \xi_{a_2} \big)
\,\,\oplus \,\,\bbC \xi_{x_3}.
\]

We can construct the finite 2-graph $\Lambda = \Lambda_G
\ast_\theta \Lambda_G$ described in $\S$~\ref{S:kgraphs} as
follows:
\[
\Lambda^{0} = \{ x_1, x_2, x_3 \},
\]
\[
\Lambda^{(1,0)} = \{ a_1, a_2 \}, \quad \Lambda^{(0,1)} = \{ b_1,
b_2\},
\]
\[
\Lambda^{(1,1)} = \{ b_2a_1\} = \{ a_2 b_1\},
\]
with range and source maps such that $x_1= s(a_1) = s(b_1)$,
$r(a_1) = x_2 = s(a_2)$, $r(b_1) = x_2 = s(b_2)$, $x_3 = r(a_2) =
r(b_2)$. The rest of $\Lambda$ consists of $\Lambda^{(2,0)} = \{
a_2a_1\}$ and $\Lambda^{(0,2)} = \{ b_2b_1\}$. Thus $\Lambda^1 =
\Lambda^{(1,0)} \cup \Lambda^{(0,1)}$ includes two `red' and two
`blue' edges and the relation $b_2a_1 = a_2b_1$ specifies all
possible commutation relations within $\Lambda$.

%\begin{figure}[h]\caption{}

%\setlength{\unitlength}{.007in}

%\begin{picture}(100,200)(300,40)

%\put(145,115){$\bullet$}

%\put(345,115){$\bullet$}

%\put(545,115){$\bullet$}

%\put(120,115){$x_1$}

%\put(345,95){$x_2$}

%\put(560,115){$x_3$}

%%%%%%%%%%%%%%%%%%%%%%%%%%%%%%%%

%\qbezier(150,120)(250,200)(350,120)

%\qbezier(150,120)(250,40)(350,120)

%\qbezier(350,120)(450,200)(550,120)

%\qbezier(350,120)(450,40)(550,120)

%%%%%%%%%%%%%%%%%%%%%%%%%%%%%%%%

%\put(250,185){$a_1$}

%\put(250,45){$b_1$}

%\put(450,185){$a_2$}

%\put(450,45){$b_2$}

%%%%%%%%%%%%%%%%%%%%%%%%%%%%%%%%

%\put(240,155){$\gg$}

%\put(440,155){$\gg$}

%\put(242,75){$>$}

%\put(442,75){$>$}

%\end{picture}

%\end{figure}

The Fock space $\H_\Lambda$ is naturally identified  with the
vertices of three disjoint downward directed graphs, with vertices
for the basis vectors $\{ \xi_{x_1},\xi_{x_2}, \xi_{x_3}\}$ at
level one, for $\{ \xi_{a_1}, \xi_{b_1}, \xi_{a_2}, \xi_{b_2}\}$
at level two and for $\{ \xi_{a_2a_1}, \xi_{b_2a_1}=\xi_{a_2b_1},
\xi_{b_2b_1}\}$ at level three. The action of $L_\lambda$,
$\lambda\in\Lambda$, is given as the appropriate downward
(partial) shift. In particular, $\fL_\Lambda$ can be identified as
a matrix algebra on a ten dimensional Hilbert space.

%\begin{figure}[h]\caption{}

%\setlength{\unitlength}{.007in}

%\begin{picture}(100,200)(300,40)

%\put(145,115){$\bullet$}

%\put(345,115){$\bullet$}

%\put(545,115){$\bullet$}

%\put(120,115){$x_1$}

%\put(345,95){$x_2$}

%\put(560,115){$x_3$}

%%%%%%%%%%%%%%%%%%%%%%%%%%%%%%%%

%\qbezier(150,120)(250,200)(350,120)

%\qbezier(150,120)(250,40)(350,120)

%\qbezier(350,120)(450,200)(550,120)

%\qbezier(350,120)(450,40)(550,120)

%%%%%%%%%%%%%%%%%%%%%%%%%%%%%%%%

%\put(250,185){$a_1$}

%\put(250,45){$b_1$}

%\put(450,185){$a_2$}

%\put(450,45){$b_2$}

%%%%%%%%%%%%%%%%%%%%%%%%%%%%%%%%

%\put(240,155){$\gg$}

%\put(440,155){$\gg$}

%\put(242,75){$>$}

%\put(442,75){$>$}

%\end{picture}

%\end{figure}

\end{eg}

\begin{eg}\label{cycleeg} ({\it Higher rank cycle algebras})
Let $C_n$ be the directed cycle graph with $n$ edges $e_i =
(x_{i+1}, x_i)$, $1\leq i \leq n$ ($i+1 \mod n$). We define
$k$-graphs $C_n^{(k)}$ which are the higher rank variants of these
graphs and identify their operator algebras as matrix function
algebras. Assume first that $k=2$. Define $C_n^{(2)}$ to be the
2-graph $\Lambda$ such that
\[
\Lambda^{(0)} = \{x_1, \ldots ,x_n\},
\]
\[
\Lambda^{(1,0)} = \{e_1, \ldots ,e_n\} \qand \Lambda^{(0,1)} =
\{f_1, \ldots ,f_n\},
\]
where $f_i$ and $e_i$ have the same sources and same ranges, and
where
\[
 f_{i+1} e_i = e_{i+1} f_i \qfor 1\leq i \leq n.
\]
In fact $\Lambda$ is the unique 2-graph arising from the
$\ast_\theta$ construction with $A=B=C_2$. The Fock space
$\H_\Lambda$ has a basis $\{\xi_\lambda\}$ which is in natural
correspondence with the vertices of $n$ disjoint graphs, each of
which is a downward directed rectangular lattice. The generators
$L_{e_i}, L_{f_i}$ can be realized as downward partial shifts,
with leftward and rightward actions, respectively. Each vertex
carries a label of the form $\xi_\lambda$ where
\[
\lambda = f_{p+q} \cdots f_{p+1} e_p \cdots e_{i+1} e_i.
\]
Thus we may identify $\H_\Lambda$ with $n$ copies of $H^2\otimes
H^2 = H^2(z,w)$, the Hardy space for the torus $\bbT^2 = \{ (z,w)
: |z|=|w|=1\}$, with its basis $\{z^p w^q : p,q \in \bbZ_+\}$.

However, there is a more useful related $n$-fold decomposition of
$\H_\Lambda$. We first illustrate this in the case $n=3$. In this
case the identification above is $\H_\Lambda \cong \bbC^3 \otimes
H^2(z,w)$ with orthonormal basis
\[
\{ g_i \otimes z^p w^q : 1\leq i \leq 3,\, p,q\in\bbZ_+ \}
\]
where $\{ g_1,g_2,g_3\}$ is an orthonormal basis for $\bbC^3$.
Consider now the decomposition $\H_\Lambda = \H_1 \oplus \H_2
\oplus \H_3$ where
\begin{align*}
\H_1 &= \spn \big\{ g_i \otimes z^{p_i}w^{q_i} : p_1 + q_1 \equiv
0, \,\, p_2 + q_2 \equiv 2, \,\, p_3+ q_3 \equiv 1 \big\}
\\
\H_2 &= \spn \big\{ g_i \otimes z^{p_i}w^{q_i} : p_2 + q_2 \equiv
0, \,\, p_3 + q_3 \equiv 2, \,\, p_1+ q_1 \equiv 1 \big\}
\\
\H_3 &= \spn \big\{ g_i \otimes z^{p_i}w^{q_i} : p_3 + q_3 \equiv
0, \,\, p_1 + q_1 \equiv 2, \,\, p_2+ q_2 \equiv 1 \big\},
\end{align*}
where each of these subspaces is closed and the prescribed
addition is modulo 3. Let us dispense with the $\bbC^3\otimes
H^2(z,w)$ identification above and identify $\H_1, \H_2, \H_3$
afresh with $H^2(\bbT^2)$ in the natural way. Now $\H_\Lambda =
H^2(\bbT^2) \oplus H^2(\bbT^2) \oplus H^2(\bbT^2)$ and we see that
the operators $L_{e_1}$, $L_{e_2}$, $L_{e_3}$, $L_{f_1}$,
$L_{f_2}$, $L_{f_3}$ are represented by the operator matrices
\[
\left[ \begin{matrix} 0&0&0\\ T_z&0&0\\ 0&0&0 \end{matrix}
\right], \quad\quad \left[ \begin{matrix} 0&0&0\\ 0&0&0\\ 0&T_z&0
\end{matrix} \right], \quad\quad \left[ \begin{matrix} 0&0&T_z\\ 0&0&0\\
0&0&0 \end{matrix} \right],
\]
\[
\left[ \begin{matrix} 0&0&0\\ T_w&0&0\\ 0&0&0 \end{matrix}
\right], \quad\quad \left[ \begin{matrix} 0&0&0\\ 0&0&0\\ 0&T_w&0
\end{matrix} \right], \quad\quad \left[ \begin{matrix} 0&0&T_w\\ 0&0&0\\
0&0&0 \end{matrix} \right],
\]
while $\alpha L_{x_1} + \beta L_{x_2} + \gamma L_{x_3}$ is
represented by
\[
\left[ \begin{matrix} \alpha I &0&0\\ 0&\beta I &0\\ 0&0&\gamma I
\end{matrix} \right].
\]
It follows readily now that $\fL_\Lambda$ is unitarily equivalent
to the matrix function algebra
\[
\left[ \begin{matrix} H^\infty_{3,0}(z,w)& &H^\infty_{3,2}(z,w)& &
H^\infty_{3,1}(z,w) \\ & &
\\
H^\infty_{3,1}(z,w)& & H^\infty_{3,0}(z,w)& &H^\infty_{3,2}(z,w)\\
& &
\\ H^\infty_{3,2}(z,w)& & H^\infty_{3,1}(z,w)& &H^\infty_{3,0}(z,w)
\end{matrix} \right]
\]
where $H^\infty_{3,i}(z,w)$ is the closed span of the basis
elements $\{z^p w^q : p+q \equiv i \mod 3 \}$, for $i=1,2$.

For the general case, $\Lambda = C_n^{(k)}$ (with $n\neq 3$,
$k\neq 2$) we have $k$ sets of morphisms/edges of total
$\delta$-degree 1, say
\[
\Lambda^{(1,0,\ldots,0)} = \{e_1^1,\ldots ,e_n^1\}, \,\,\ldots
,\,\, \Lambda^{(0,\ldots,0,1)} = \{e_1^k,\ldots ,e_n^k\}
\]
and all other morphisms arise from compositions, subject only to
identifications through the relations
\[
 e_{i+1}^s e_i^r = e_{i+1}^r e_i^s
\]
for all $i$ ($i+1 \mod n$), and all $1 \leq r\neq s \leq k$. We
identify the subspace of $\H_\Lambda$ which is spanned by
$\{\lambda : r(\lambda) = x_i \}$ with $H^2(\bbT^k)$. Then
$\H_\Lambda$ is isomorphic to $H^2(\bbT^k) \oplus \ldots \oplus
H^2 (\bbT^k)$ and we identify, as before, $\fL_\Lambda$ with the
matrix function algebra
\[
\left[ \begin{matrix}  H^\infty_{n,0} (\bbT^k)&H^\infty_{n,n-1} (\bbT^k)& \cdots &H^\infty_{n,1} (\bbT^k)\\
H^\infty_{n,1} (\bbT^k)&H^\infty_{n,0} (\bbT^k)&   & \vdots \\
\vdots & & \ddots & \\
H^\infty_{n,n-1} (\bbT^k)& \cdots & & H^\infty_{n,0} (\bbT^k)
\end{matrix} \right],
\]
where $H^\infty_{n,i} (\bbT^k)$ is the weak-$\ast$ closed subspace
of $H^\infty (\bbT^k)$ spanned by the monomials $z_1^{i_1} \cdots
z_k^{i_k}$ with $i_1 + \ldots + i_k \equiv i \mod n$.

In view of the matrix function identifications in Alaimia and
Peters \cite{AlPe}, it is now possible to see that for the higher
rank cycle graph $C_n^{(k)}$ its higher rank semigroupoid algebra
is equal to the higher rank $\sigma$-weakly closed semicrossed
product $\bbC^n \times_\alpha^\sigma \bbZ_+^k$ where the action
$\alpha : \bbZ_+^k \rightarrow \Aut(\bbC^n) $ is given by
$\alpha(m_1,\ldots,m_k) = \sigma^{m_1+\ldots + m_k}$ where
$\sigma$ is the cyclic shift.
\end{eg}

\section{The Algebras $\A_{\underline{n},\theta}$ and
$\fL_{\underline{n},\theta}$}\label{S:singlevertex} %%%%%%%%%%%%%%%%%%%%%%%%%%%

We now consider single vertex $k$-graphs and their nonselfadjoint
operator algebras $\A_{\underline{n},\theta}$ (norm closed) and
$\fL_{\underline{n},\theta}$ ($\wot$-closed). First we determine
the codimension one invariant subspaces of
$\fL_{\underline{n},\theta}$ and identify the natural connection
with the Gelfand space of the quotient of the algebra
$\A_{\underline{n},\theta}$ by its norm-closed commutator ideal.
As we shall see this Gelfand space is biholomorphically equivalent
to a subspace of the direct product $\bbB_{n_1} \times \ldots
\times \bbB_{n_k}$ determined by a complex algebraic variety
associated with the relations latent in the $k$-graph.

Let us now specify a general single vertex $k$-graph explicitly in
terms of edge generators and commutation relations. We write
$\Lambda_{\underline{n},\theta}$ for such a $k$-graph, where
$\underline{n} = (n_1,\ldots ,n_k)$, where $n_i$ is the number of
edges $e$ of degree $d(e) = \delta_i = (0,\ldots,0,1,0,\ldots,0)$,
and where $\theta$ denotes a set $\{\theta_{i,j}:1\leq i<j\leq
k\}$ of permutations which determine the relations
\[
e^{(i)} e^{(j)} = \big(
\theta_{i,j}(e^{(i)}e^{(j)})\big)^{\rm{op}},
\]
where $d(e^{(i)}) = \delta_i$, and where $(ef)^{\rm{op}}$ denotes
the opposite product $fe$. Thus $\theta_{i,j}$ is a permutation of
the $n_i n_j$ products $e^{(i)} e^{(j)}$, which, when necessary,
we enumerate in the natural order
\[
e_1^{(i)} e_1^{(j)},e_1^{(i)} e_2^{(j)}, \ldots e_1^{(i)}
e_{n_j}^{(j)}, e_2^{(i)} e_1^{(j)}, \ldots \ldots ,e_{n_i}^{(i)}
e_{n_j}^{(j)}.
\]
Here we have labelled the edges of degree $\delta_i$ as
$e_1^{(i)}, e_2^{(i)},\ldots, e_{n_i}^{(i)}$.

Consider now a path $\lambda$ in $\Lambda =
\Lambda_{\underline{n},\theta}$ with unique factorization $\lambda
= \lambda_1 \lambda_2 \cdots \lambda_k$ where each $\lambda_j$ is
a free word in the loop edges of degree $\delta_j$. For a point
$\alpha^{(j)} \in \bbC^{n_j}$ define the corresponding word
$\lambda_j(\alpha^{(j)})$, corresponding to evaluation in $\bbC$
of the word $\lambda_j$, by letterwise substitution, at
$\alpha^{(j)} = (\alpha_1^{(j)},\ldots,\alpha_{n_j}^{(j)})$ in
$\bbC^{n_j}$. Finally, for $\alpha =
(\alpha^{(1)},\alpha^{(2)},\ldots,\alpha^{(k)})$ in
$\bbC^{n_1}\times \ldots \times \bbC^{n_k}$ define
\[
\lambda(\alpha) = \lambda_1(\alpha^{(1)})\lambda_2(\alpha^{(2)})
\ldots \lambda_k(\alpha^{(k)}).
\]
Thus we only evaluate the general path $\lambda$ if it is
expressed in its uniquely factored form. If $\alpha$ lies in the
open ball product $\bbB_{n_1}^0 \times \ldots \times \bbB_{n_k}^0$
then we may define the unit vector $\nu_\alpha = \omega_\alpha /
||\omega_\alpha||_2$ in the Fock space $\H_\Lambda$, where
$\omega_\alpha = \sum_{\lambda\in\Lambda} \lambda(\alpha)
\xi_\lambda$. Indeed,
\begin{eqnarray*}
||\omega_\alpha||_2^2 &=& \sum_{\lambda\in\Lambda}
|\lambda(\alpha)|^2 \\
&=& \sum_{\lambda_1\in\bbF^+_{n_1}} \ldots
\sum_{\lambda_k\in\bbF^+_{n_k}} |\lambda_1(\alpha^{(1)})|^2 \ldots
|\lambda_k(\alpha^{(k)})|^2 \\
&=& \prod_{i=1}^k \big( 1 - ||\alpha^{(i)}||_2^2\big)^{-1}.
\end{eqnarray*}

Suppose first that the commutation relations given by $\theta = \{
\theta_{i,j} : 1 \leq i < j \leq k\}$ are the commuting relations
arising when each $\theta_{i,j}$ is the identity permutation of
$d^{-1}(\delta_i)d^{-1}(\delta_j)$. In particular, for each
generating edge $e = e_j^{(i)}$ with degree $\delta_i$ we have
\[
e\lambda = e \lambda_1\lambda_2 \cdots \lambda_k = \lambda_1
\cdots \lambda_{i-1} (e\lambda_i)\lambda_{i+1} \cdots \lambda_k
\]
and so $(e\lambda)(\alpha) = \alpha_j^{(i)}\lambda(\alpha)$. We
may now deduce that $L_e^* \omega_\alpha = \alpha_j^{(i)}
\omega_\alpha$. Indeed, for all $\lambda\in \Lambda$,
\begin{eqnarray*}
\big< L_e^* \omega_\alpha , \xi_\lambda \big> &=& \big<
\omega_\alpha , \xi_{e\lambda} \big> = (e\lambda)(\alpha) \\
&=& \alpha_j^{(i)} \lambda(\alpha) = \alpha_j^{(i)} \big<
\omega_\alpha , \xi_\lambda \big> = \big< \alpha_j^{(i)}
\omega_\alpha , \xi_\lambda \big>.
\end{eqnarray*}
Thus, with $N=n_1+\ldots + n_k$, we have shown that each $N$-tuple
in the product $\bbB_{n_1}^0 \times \ldots \times \bbB_{n_k}^0$
(of open unit balls) is a joint eigenvalue for the $N$-tuple
$\big\{ L_{e_1^{(1)}}^*, \ldots , L^*_{e_{n_k}^{(k)}}\big\}$ with
eigenvector $\omega_\alpha$, and that $\{\omega_\alpha\}^\perp$ is
therefore a codimension one subspace in $\Lat \fL_\Lambda$. Here,
as we have already noted in (\ref{1steg}),
$\fL_{\underline{n},\theta}$ is naturally identifiable with the
spatial tensor product $\fL_{n_1} \otimes \ldots \otimes
\fL_{n_k}$.

Suppose now that $\theta = \{\theta_{i,j} : 1 \leq i < j \leq k\}$
is a general family of permutations. For $1\leq i \leq k$ let
$z_{i,1}, \ldots , z_{i,n_i}$ be the coordinate variables for
$\bbC^{n_i}$ so that there is a natural bijective correspondence
$e_k^{(i)} \rightarrow z_{i,k}$ between edges and coordinate
variables. We define $V_\theta \subseteq \bbC^{n_1} \times \ldots
\times \bbC^{n_k}$ to be the complex algebraic variety determined
by  the equation set
\[
\big\{ z_{i,p}z_{j,q} - \hat{\theta}_{i,j}(z_{i,p}z_{j,q}) : 1\leq
p \leq n_i,\, 1\leq q \leq n_j, \,1\leq i<j\leq k \big\}
\]
where $\hat{\theta}_{i,j}$ is the permutation induced by
$\theta_{i,j}$ and the bijective correspondence.

We now obtain the following identification of the eigenvalues for
the adjoint algebra of $\fL_{\underline{n},\theta}$ and the
Gelfand space of the norm closed algebra
$\A_{\underline{n},\theta}$. We write $\bbB_{\underline{n}}^0$ for
the product of open unit balls $\bbB_{n_1}^0 \times \ldots \times
\bbB_{n_k}^0$.

\begin{thm}\label{gelfand}
(i) Each invariant subspace of $\fL_{\underline{n},\theta}$ of
codimension one has the form $\{\omega_\alpha\}^\perp$ for some
$\alpha $ in $\bbB_{\underline{n}}^0 \cap V_\theta$.

{\noindent}(ii) The character space
$\M(\A_{\underline{n},\theta})$ is biholomorphically isomorphic to
$\bbB_{\underline{n}}\cap V_\theta$ under the map $\varphi$ given
by
\[
\varphi (\rho) = \big(\rho(L_{e_1^{(1)}}),
\ldots,\rho(L_{e^{(k)}_{n_k}})\big), \qfor
\rho\in\M(\A_{\underline{n},\theta}).
\]
\end{thm}

\Prf To see $(ii)$ first let $\rho \in
\M(\A_{\underline{n},\theta})$. Since $\rho$ is a multiplicative
linear functional it is completely contractive. Thus, since the
row operator $R_i = [L_{e_1^{(i)}}\,\, L_{e_2^{(i)}}\, \cdots
\,L_{e_{n_i}^{(i)}}]$ satisfies $R_i R_i^*\leq I$ it follows that
the scalar row matrix $[\rho(L_{e_1^{(i)}})\,\,
\rho(L_{e_2^{(i)}}) \,\cdots\, \rho(L_{e_{n_i}^{(i)}})]$ is a
contraction, and hence that $\alpha^{(i)} =
(\rho(L_{e_k^{(i)}}))_{k=1}^{n_i}$ lies in $\bbB_{n_i}$. Let
$\alpha = (\alpha^{(1)},\ldots ,\alpha^{(k)})$ be the point in
$\bbB_{\underline{n}}$ which derives in this way from $\rho$. We
have thus shown that the map $\varphi$ maps
$\M(\A_{\underline{n},\theta})$ into the product ball
$\bbB_{n_i}$. In view of the relations $e_p^{(i)} e_q^{(j)} =
(\theta_{i,j}(e_p^{(i)} e_q^{(j)}))^{\rm{op}}$, we have $e_p^{(i)}
e_q^{(j)} = e_s^{(j)} e_r^{(i)}$ for some $r,s$ depending on
$p,q$, thus
\begin{eqnarray*}
\alpha_p^{(i)} \alpha_q^{(j)} &=& \rho (L_{e_p^{(i)}}
L_{e_q^{(j)}}) =\rho (L_{e_p^{(i)}} L_{e_q^{(j)}}) \\ &=& \rho
(L_{e_p^{(i)} e_q^{(j)}}) = \rho (L_{e_s^{(j)} e_r^{(i)}}) \\ &=&
\rho (L_{e_s^{(j)}}) \rho (L_{e_r^{(i)}}) = \alpha_s^{(j)}
\alpha_r^{(i)}
\end{eqnarray*}
and so in particular the polynomial
\[
z_{i,p} z_{j,q} - \hat{\theta}_{i,j} (z_{i,p}z_{j,q}) = z_{i,p}
z_{j,q} - z_{j,s} z_{i,r}
\]
vanishes on $\alpha$. This is true for all appropriate $p,q,i,j$
and so $\alpha \in \bbB_{\underline{n}} \cap V_\theta$.

On the other hand suppose that $\alpha\in V_\theta$. Then it
follows that for $e=e_j^{(i)}$ we have $(e\lambda) (\alpha) =
\alpha_j^{(i)} \lambda(\alpha)$, and indeed, for any (unfactored)
path $\lambda$ in the edges of the $k$-graph the substitutional
evaluation of $\lambda$ at $\alpha$ coincides with the evaluation
of the factored form of $\lambda$ at $\alpha$, which we denote
$\lambda(\alpha)$. If, in addition, $\alpha\in
\bbB_{\underline{n}}^0 \cap V_\theta$ then our earlier calculation
shows that the vector $\omega_\alpha$ is an eigenvector for the
joint eigenvalue $\alpha$ for the $N$-tuple  $(L^*_{e_1^{(1)}},
\ldots ,L^*_{e_{n_k}^{(k)}})$. It follows readily that the unit
vector $\nu_{\ol{\alpha}}$ defines a vector functional
\[
\rho(A) = \left< A\nu_{\ol{\alpha}}, \nu_{\ol{\alpha}}\right>
\]
which defines a character $\rho$ in
$\M(\A_{\underline{n},\theta})$ with $\varphi(\rho) = \alpha$.
Since $\M(\A_{\underline{n},\theta})$ is a compact Hausdorff
space, and since we have shown that the range of $\varphi$
contains $\bbB_{\underline{n}}^0\cap V_\theta$ and is contained in
$\bbB_{\underline{n}}^0\cap V_\theta$, it follows that the range
of $\varphi$ is precisely    $\bbB_{\underline{n}}^0\cap
V_\theta$. From this, part $(ii)$ of the theorem now follows.

To see $(i)$, suppose now that $\nu$ is a unit vector such that
$\{\nu\}^\perp$ is invariant for $\fL_{\underline{n},\theta}$, and
hence that $\nu$ is a joint eigenvector for the
$|\underline{n}|$-tuple $(L_{e_1^{(1)}}^*,\ldots ,
L_{e_{n_k}^{(k)}}^*)$ with corresponding eigenvalue $\beta = (
\beta^{(1)}, \ldots , \beta^{(k)})$. Since the column operators
$R_i^*$, $1\leq i \leq k$, are contractions it follows that
$\beta^{(i)}\in\bbB_{n_i}$. Since $\nu$ is a joint eigenvector it
follows that the map $L_e \mapsto \langle L_e \nu,\nu \rangle$
extends to a multiplicative linear functional on
$\A_{\underline{n},\theta}$ and so from the calculation above
$\ol{\beta}$, and hence $\beta$, lies in $\bbB_n\cap V_\theta$.

Suppose now that $\nu = \sum_{\lambda\in\Lambda} b_\lambda
\xi_\lambda$. Then
\[
b_\lambda = \langle \nu, \xi_\lambda \rangle = \langle L_\lambda^*
\nu, \xi_x \rangle \\ = \lambda(\beta) \langle \nu, \xi_x \rangle
= \lambda(\beta) b_x,
\]
where, as before, $\lambda(\beta)$ represents the factored form
substitution of the word $\lambda$ in $\Lambda$. Our earlier
calculation of the norm of $\nu_\alpha$ applies here and the
finiteness of $\sum |\lambda(\beta)|^2$ implies that $\beta\in
\bbB_{n_1}^0 \cap V_\theta$. \bx

\begin{rem}
As an illustration of the theorem let $\underline{n} = (n_1,n_2)$
and let $z_1,\ldots,z_{n_1}$ and $w_1,\ldots,w_{n_2}$ be
coordinate variables for $\bbC^{n_1}$ and $\bbC^{n_2}$. If
$\theta$ is a simple cyclic permutation of all the $n_1n_2$ pairs
$\{z_iw_j\}$ then $V_\theta$ is the variety $V_\theta = \{ (z,w):
z_1=\ldots =z_{n_1}, \, w_1=\ldots = w_{n_2}\}$ and so
$\bbB_{\underline{n}}\cap V_\theta$ is the direct product of two
discs with radii $n_1^{- 1 / 2}$, $n_2^{- 1 / 2}$. This is the
minimal such subset of $\bbB_{\underline{n}}$ associated with a
permutation and it is not hard to see that many other
permutations, including products of cycles, also lead to this
minimal case. In general it can be shown that the Gelfand space
$\bbB_{\underline{n}}\cap V_\theta$ does not determine the
operator algebra $\A_{\underline{n},\theta}$ up to isometric
isomorphism.
\end{rem}

%%%%%%%%%%%%%%%%%%%%%%%%%%%%%%%%%%%%%%%%%%%%%%%%%%%
\section{Reflexivity}\label{S:reflexivity}
%%%%%%%%%%%%%%%%%%%%%%%%%%%%%%%%%%%%%%%%%%%%%%%%%%%

Recall that a ($\wot$-closed) operator algebra $\fA$ is {\it
reflexive} if $\fA$ coincides with the algebra of operators which
leave every subspace in the invariant subspace lattice for $\A$
invariant, $\fA = \Alg \Lat \fA$. On the other hand, a measure of
the distance to an operator algebra $\fA$ is given by $
\beta_{\fA} (X) = \sup_{L\in\Lat\fA} ||P_L^\perp X P_L||, $ where
$P_L$ is the projection onto the subspace $L$. Clearly $\beta_\fA
(X) \leq \dist (X,\fA)$ and $\fA$ is said to be {\it
hyper-reflexive} if there is a constant $C$ such that
$\dist(X,\fA) \leq C \beta_\fA (X) $ for all $X$. We begin by
identifying a new class of hyper-reflexive algebras.

As a generalization of terminology from \cite{KP1,KP2}, we define
the `double pure cycle property' for a higher rank graph. Firstly,
a {\it pure cycle} is one composed of edges of the same degree (a
monochromatic cycle), and, secondly, $\Lambda$ has the {\it double
pure cycle property} if for every $v\in\Lambda^0$ there is a path
$\lambda\in\Lambda$ with $s(\lambda)=v$ and $r(\lambda)=w$ such
that $w$ lies on a {\it double pure cycle} in the sense that there
is a pair of distinct pure cycles $\lambda_i = w\lambda_i w$,
$i=1,2$, of the same colour, neither of which may be written as a
product of cycles.

\begin{lem}\label{dclemma}
If $\Lambda$ satisfies the double pure cycle property then
$\fL_\Lambda$ contains a pair of isometries with mutually
orthogonal ranges.
\end{lem}

\Prf We may construct isometries $U,V\in\fL_\Lambda$ with $U^*V=0$
in a direct manner as follows. Let $\lambda_1\neq \lambda_2$ be a
double pure cycle with $s(\lambda_i)=r(\lambda_i)=v$ for some
$v\in\Lambda^0$.
%Without loss of generality we will
We may assume that
for all $w\in\Lambda^0$ there is a $\lambda_w\in\Lambda$ such that
$s(\lambda_w)=w$ and $r(\lambda_w)=v$.
%In the general case a
%`splicing together' construction may be used to build the desired
%isometries via pairs of partial isometries associated with each
%vertex.
since the general case follows easily from this special case.
By hypothesis and from the factorization property, for $k\geq 1$
the paths $\lambda_1^k\lambda_2$ are cycles over $v$ and the
partial isometries $L_{\lambda_1^k\lambda_2}$, $k\geq 1$, have
mutually orthogonal ranges with initial projection $L_v$. Let
$w\mapsto\{k_a^w,k_b^w\}$ be a one-to-two map from $\Lambda^0$ to
the positive integers $\bbN$. As the desired isometries we may
define
\[
U = \sum_{w\in\Lambda^0} L_{\lambda_1}^{k_a^w} L_{\lambda_2}
L_{\lambda_w} \quad\qand\quad V = \sum_{w\in\Lambda^0}
L_{\lambda_1}^{k_b^w} L_{\lambda_2} L_{\lambda_w},
\]
the sums converging $\sot$ when $\Lambda^0$ is infinite. \bx

\begin{thm}\label{hrthm}
If $\Lambda^t$ satisfies the double pure cycle property then
$\fL_\Lambda$ is hyper-reflexive with distance constant at most 3.
\end{thm}

\Prf As $\fL_{\Lambda^t}$ is unitarily equivalent to $\fR_\Lambda
= \fL_\Lambda^\prime$, the previous lemma shows that
$\fL_\Lambda^\prime$ contains a pair of isometries with mutually
orthogonal ranges. Thus the result follows as a direct application
of Bercovici's hyper-reflexivity Theorem \cite{Berc}. \bx

As an immediate consequence we obtain the following.

\begin{cor}\label{svhr}
Let $\underline{n}=(n_1,\ldots,n_k)\in\bbN^k$ and suppose that
$n_j\geq 2$ for some $j$. Then $\fL_{\underline{n},\theta}$ is
hyper-reflexive for all choices of $\theta$.
\end{cor}

Note that the single vertex algebras $\fL_\Lambda$ which do not
satisfy the hypothesis of Corollary~\ref{svhr} are each unitarily
equivalent to $\bbC$ or $H^\infty(\bbT^k) \cong
(H^\infty)^{\otimes k}$ for some $k\geq 1$, and these algebras are
known to be reflexive \cite{Sarason}. (We note that the problem of
hyper-reflexivity for $H^\infty(\bbT^k)$, $k\geq 2$, appears to
remain unresolved at present.) Reflexivity of these algebras is
well-known, but for the interested reader we mention that this
fact may be deduced from the first part of the proof of
Theorem~\ref{lnreflexive}. Thus, as hyper-reflexivity subsumes
reflexivity, it follows now
that every single vertex algebra
$\fL_\Lambda$ is reflexive. We use this in the proof below.

We shall prove reflexivity for $\fL_\Lambda$ up to a mild graph
constraint. We shall say $v\in\Lambda^0$ is a {\it radiating
vertex} when $\lambda\in\Lambda^1$ with $r(\lambda)=v$ implies
that $s(\lambda)=v$. Such a vertex is {\it multiplicity one} if
there is at most one loop edge at $v$ of each colour. Further
we say
that a radiating vertex $v\in\Lambda^0$ is {\it relational} if
there are loop edges $\mu\neq\mu'$ at $v$ and paths
$\lambda,\lambda'\in\Lambda$ with $s(\lambda)=v=s(\lambda')$ that
immediately leave $v$ such that $\lambda\mu = \lambda'\mu'$.

\begin{thm}\label{reflexivethm}
Let $\Lambda$ be a higher rank graph with no multiplicity one
relational radiating vertices. Then $\fL_{\Lambda}$ is reflexive.
\end{thm}

\Prf Let $A\in \Alg \Lat \fL_\Lambda$. We shall show that if $x\in
d^{-1}(0)$ and $A = AL_x$ then $A\in \fL_\Lambda$. Since every
operator $B$ on $\H_\Lambda$ is the weak operator topology limit
of the sums $\sum_{i\geq 1} BL_{x_i}$, where $x_1, x_2, \ldots$ is
an enumeration of $d^{-1}(0)$, the proof will be complete.

Given $\mu\in\Lambda$ with $r(\mu)=x$ we have
\[
A\xi_\mu = \sum_{s(\lambda)=x} \alpha_\lambda^\mu \xi_{\lambda\mu}
\]
for some choice of scalars $\alpha_\lambda^\mu$. This follows
since the subspace $\M$ spanned by $\{ \xi_{\lambda\mu} :
\lambda\in\Lambda\}$ belongs to $\Lat\fL_\Lambda$. Note that
$A\xi_\mu=0$ if  $r(\mu)\neq x$. We shall show that for all paths
$\mu,\nu$ with $r(\mu)=x=r(\nu)$ we have $\alpha_\lambda^\mu =
\alpha_\lambda^\nu$ for all paths $\lambda$ with $s(\lambda)=x$.
If this is the case then for all $\lambda'\in\Lambda$ with
$r(\lambda') = s(\mu)$
\[
R_{\lambda'} A \xi_\mu = \sum_\lambda \alpha_\lambda^\mu
\xi_{\lambda\mu\lambda'} = \sum_\lambda
\alpha_\lambda^{\mu\lambda'} \xi_{\lambda\mu\lambda'} = A
\xi_{\mu\lambda'} = A R_{\lambda'} \xi_\mu,
\]
and $R_{\lambda'}A\xi_\mu = 0 = AR_{\lambda'}\xi_\mu$ when
$r(\lambda')\neq s(\mu)$.  Hence $A\in\fR_\Lambda^\prime$ and so
$A\in\fL_\Lambda$, as desired.

We consider two cases. Suppose first that there is a path $\nu$
with $\nu=x\nu y$ and $y\neq x$ ($x,y\in d^{-1}(0)$). Then the
range $\N$ of $R_x + R_\nu$ is spanned by the set of vectors
$\{\xi_{\lambda x} + \xi_{\lambda \nu} : s(\lambda)=x\}$ and the
vectors in this set are pairwise orthogonal. Since
$\N\in\Lat\fL_\Lambda$ it follows that
\[
A(\xi_x + \xi_\nu) = A (R_x + R_\nu) \xi_x = \sum_{s(\lambda)=x}
\gamma_\lambda (\xi_{\lambda x} + \xi_{\lambda \nu})
\]
for some choice of scalars $\gamma_\lambda$. But $A\xi_x$ and
$A\xi_\nu$ are given in terms of the coefficients
$\alpha_\lambda^x$, $\alpha_\lambda^\nu$ respectively and so
$\alpha_\lambda^x = \gamma_\lambda = \alpha_\lambda^\nu$ since
$s(\nu)\neq x$. Precisely the same argument holds if we replace
$x$ by a path $\mu$ with $\mu=x\mu x$, and so we obtain
$\alpha_\lambda^\mu = \alpha_\lambda^x = \alpha_\lambda^\nu$ for
all $\lambda$ for such a path at $x$. It follows that
$\alpha_\lambda^{\mu'} = \alpha_\lambda^x$ for all $\lambda$ and
for all paths $\mu'$ which terminate at $x$, as desired.

If there is no such path $\nu =x\nu y$ with $y\neq x$ then we are
in the second case in which $\mu=x\mu x$ whenever $r(\mu)=x$.
Plainly this entails that there is a single vertex induced
sub-$k$-graph $\Gamma$ of $\Lambda$ such that $L_x
\fL_\Lambda|_{L_x \H_\Lambda}$ is unitarily equivalent to
$\fL_\Gamma$. By the previous discussion $\fL_\Gamma$ is
reflexive, and so we do at least have $\alpha_\lambda^\mu =
\alpha_\lambda^x$, for all $\lambda =x\lambda x$, when $\lambda =x
\lambda $ ($=x\lambda x$). We shall now show that this equality
also holds for paths $\lambda'$ with $\lambda' = y\lambda' x$,
$y\neq x$, and this will complete the proof.

If the loop edges at $x$ include a double pure loop then we may
argue as above and use the Bercovici Theorem to deduce this
equality for all $\lambda^\prime$. Further, the equality trivially
holds when there are no loops at $x$. Thus we may reduce to the
case that $x$ is a multiplicity one vertex.

Suppose first that $\lambda'$ is not of the form $\lambda_1 h$
with $h\in \Gamma$ and with $d(h)\neq 0$. Consider the restriction
operator $A_{\lambda'} = L_{\lambda'}^* A|_{\H_\Gamma}$. We show
that $A_{\lambda'}$ lies in $\fL_\Gamma$. To this end let $\M\in
\Lat\fL_\Gamma$ and define $\td{\M} = \bigvee_{s(\lambda)=x}
L_\lambda \M$ in $\Lat \fL_\Gamma$. Then
\[
A_{\lambda'} \M = L_{\lambda'}^* A \M \subseteq L_{\lambda'}^*
\td{\M}.
\]
In view of the assumption on $\lambda'$, and the constraint on
$\Lambda$ in the hypothesis, we have that $L_{\lambda'}^*
L_{\lambda}$ is non-zero for $s(\lambda)=x$ only if $\lambda =
\lambda'h$ with $h\in\Gamma$.
%Indeed, the range of $L_{\lambda'}$ is spanned by
%$\{\xi_{\lambda'h}:h\in\Gamma\}$, and similarly for $L_w$, so the
%assumption on $w'$ implies an inner product
%$(L_w\xi_h,L_{w'}\xi_{h'})$ can only be non-zero if $w=w'h$ in
%this case.
Thus
\[
L_{\lambda'}^*\td{\M} \subseteq \bigvee_{h=xhx} L_h \M \subseteq
\M.
\]
Since we have shown that $A_{\lambda'}$ belongs to $\Alg
\Lat\fL_\Gamma$, it follows that $A_{\lambda'}$ is in $\fL_\Gamma$
and hence that there exist scalars $\alpha_h$ such that
$A_{\lambda'} \sim \sum_{h=xhx} \alpha_h L_h$.

Thus if $\mu$ is a path in $\Gamma$ then
\begin{eqnarray*}
\sum_{h\in\Gamma} \alpha_h \xi_{\lambda'h\mu} &=& L_{\lambda'}
\big( A_{\lambda'}
\xi_\mu\big) \\
&=& L_{\lambda'} L_{\lambda'}^* A \xi_\mu \\
&=& L_{\lambda'} L_{\lambda'}^* \Big( \sum_{s(\lambda)=x} \alpha_\lambda^\mu \xi_{\lambda\mu} \Big) \\
&=& \sum_{h\in\Gamma} \alpha^\mu_{\lambda'h} \xi_{\lambda'h\mu}.
\end{eqnarray*}
Therefore, $\alpha^\mu_{\lambda'h} = \alpha_h$ for all
$h,\mu\in\Gamma$ and $\alpha_{\lambda'h}^\mu =
\alpha_{\lambda'h}^\nu$ for all $\mu,\nu\in\Gamma$. As we are in
the second case it follows that $\alpha_\lambda^\mu =
\alpha_\lambda^\nu$ for all $\mu,\nu$ with $r(\mu)=x=r(\nu)$ and
$\lambda$ with $s(\lambda)=x$, as desired. \bx

%%%%%%%%%%%%%%%%%%%%%%%%%%%%%%%%%%%%%%%%%%%%%%%%%%%%%

The following proof of reflexivity for the free semigroup algebras
$\fL_n$ is considerably more elementary than other proofs in the
literature \cite{ArPop,DP2} and gives a more direct generalization
of Sarason's approach for $\fL_1 = H^\infty$ \cite{Sarason}. We
include it for interest's sake.

\begin{thm}\label{lnreflexive}
$\fL_n$ is reflexive.
\end{thm}

\Prf Let $A\in \Alg \Lat \fL_n$. For $\alpha =
(\alpha_1,\ldots,\alpha_n)$ in $\bbB_n$ we may define eigenvectors
for $\fL_n^*$ by $\nu_\alpha = \sum_{w\in\bbF^+_n} \ol{w(\alpha)}
\xi_w$. Then $\{ \nu_\alpha \}^\perp$ is invariant for $\fL_n$ and
so $A^* v_\alpha = \ol{\lambda}_\alpha v_\alpha$ for some scalar
$\lambda_\alpha$ in $\bbC$.

Let $A\xi_x = \sum_w a_w \xi_w$. Then
\begin{eqnarray*}
\lambda_\alpha \big< \xi_x, v_\alpha \big> &=& \big< A \xi_x
,v_\alpha \big>
= \sum_{w} a_w \big< \xi_w , v_\alpha \big> \\
&=& \sum_{w} a_w \big< \xi_x , L_w^* v_\alpha \big>  = \sum_{w}
a_w w(\alpha) \big< \xi_x , v_\alpha \big>,
\end{eqnarray*}
and so $\lambda_\alpha = \sum_w a_w w(\alpha)$.

Note that for each $v\in \bbF^+_n$ the subspace $\M_v$ spanned by
$\{\xi_{wv}:w\in \bbF_n^+\}$ belongs to $\Lat \fL_n$ and so for
some scalars $b_w^v$, $w\in\bbF_n^+$,
\[
A\xi_v = \sum_w b_w^v \xi_{wv}.
\]
We show that $b_w = a_w$ for all $w$. This will complete the proof
since we obtain for any choice of $v$
\[
R_v A \xi_x = \sum_w a_w \xi_{wv} = \sum_w b_w^v \xi_{wv} = A
\xi_v  = A R_v \xi_x.
\]
As $A\in\Alg\Lat\fL_n$ was arbitrary, it follows that $R_v AL_w
=AL_wR_v = AR_vL_w$ and hence $R_vA\xi_w = R_vAL_w\xi_x =
AR_v\xi_w$ for all $w$. Thus $R_v A = A R_v$ for all $v$, and so
$A$ belongs to $\fR_n^\prime = \fL_n$.

Observe that
\begin{eqnarray*}
\big< A \xi_v , v_\alpha \big> &=& \sum_w b_w^v \big< L_{wv}
\xi_x,
v_\alpha \big> \\
&=& \sum_w b_w^v w(\alpha) v(\alpha) \big< \xi_x , v_\alpha \big>
.
\end{eqnarray*}
Further,
\[
\big< A \xi_v , v_\alpha \big> = \big<  \xi_v , A^* v_\alpha \big>
= \lambda_\alpha v(\alpha) \big<  \xi_x , v_\alpha \big>,
\]
and hence
\[
\lambda_\alpha = \sum_w b_w^v w(\alpha).
\]
Now
\[
\sum_w b_w^v w(\alpha) = \sum_w a_w w(\alpha)
\]
for all $||\alpha ||_2 < 1$. (Observe that in the case of
$H^\infty$ or $H^\infty(\bbT^k)$ this fact finishes the proof.)

In particular, we have $a_x = b_x^v$ for all choices of $v$. Thus,
by replacing $A$ with $A - a_xI$ we may assume that $a_x = b_x^v =
0$ for all $v$. Suppose now that $k\geq 1$ is minimal such that
$a_w = 0 = b_w^v$ for all $|w|< k$ (where $|w|= \delta(w)$ is word
length) and all $v$. We claim that $L_{w'}^* A$ belongs to $\Alg
\Lat \fL_n$ for all $|w'|=k$. If this holds then observe
\begin{eqnarray*}
L_{w'}^* A \xi_x &=& L_{w'}^* \sum_w a_w \xi_w = \sum_u a_{w'u}
\xi_u \\
L_{w'}^* A \xi_v &=& L_{w'}^* \sum_w b_w^v \xi_{wv} = \sum_u
b_{w'u}^v \xi_{uv},
\end{eqnarray*}
and hence by the above argument $a_{w'} = b_{w'}^v$ for all $v$.
Thus we finish the proof by verifying the claim.

Let $\M$ belong to $\Lat \fL_n$. Without loss of generality assume
$\M$ is cyclic. Then an orthonormal basis for $\M$ is given by $\{
L_w \eta : w\in \bbF_n^+\}$ where $\eta$ is some unit  vector.
(This is follows from part of the Beurling Theorem for $\fL_n$
\cite{DP2,Pop1}.) As $A\M \subseteq \M$ we have $A\eta = \sum_w
c_w L_w \eta$ for some scalars $c_w$. By our assumption on the
scalars $a_w, b_w^v$, a `graded Fock space' type argument can be
used to show that $c_w = 0$ for all $|w|<k$. Thus
\[
L_{w'}^* A \eta = \sum_{|w|\geq k} c_w L_{w'}^* L_w \eta  =
\sum_{u} c_{w'u}  L_u \eta \in \M.
\]
More generally, for arbitrary $u$ we have $AL_u\eta =
\sum_{|w|\geq k + |u|} c_w L_w \eta$ and similarly $L_{w'}^* A L_u
\eta \in \M$. Thus $L_{w'}^* A \M \subseteq \M$ for all $|w'|=k$
and all cyclic subspaces $\M\in \Lat\fL_n$ and it follows that
$L_{w'}^*A $ belongs to $\Alg\Lat \fL_n$, as claimed. \bx

\section{Semisimplicity}\label{S:ss}
%%%%%%%%%%%%%%%%%%%%%%%%%%%%%%%%%%%%%%%%%%%%%%%%%%%%%%%%%%%%%%%%%%%%%%%%%%%%%%%

We say that an edge $\lambda\in\Lambda^1$ `lies on a cycle' when
there is a cycle $\mu\in\Lambda$, $s(\mu)=r(\mu)$, that includes
$\lambda$ in at least one of its factorizations as a product of
edges.  Let $\nc(\Lambda)$ be the edges in $\Lambda^1$ that do not
lie on a cycle.  We show that the Jacobson radical of
$\fL_\Lambda$, $\rad(\fL_\Lambda)$, is determined by the operators
$L_\lambda$, $\lambda\in\nc(\Lambda)$, and in the finite vertex
case we obtain a complete description of $\rad(\fL_\Lambda)$.
Recall that $\rad(\fL_\Lambda)$ is the largest quasinilpotent
ideal of $\fL_\Lambda$ and that $\fL_\Lambda$ is semisimple if and
only if $\rad(\fL_\Lambda)$ is the zero ideal.

We begin with a combinatorial lemma that shows extremal paths, in
the sense of the distance measure, have special factorization
properties. Given $\lambda,\mu\in\Lambda$ with
$d(\lambda)=(m_1,\ldots ,m_k)$ and $d(\mu)=(n_1,\ldots,n_k)$ we
write $\lambda\geq \mu$ when the corresponding lexicographic
ordering on these $k$-tuples is satisfied in $\bbN^k$.

\begin{lem}\label{comblemma}
Let $\Gamma$ be a nonempty subset of $\Lambda$ such that
$\delta(\lambda_1)=\delta(\lambda_2)$ for all
$\lambda_1,\lambda_2\in\Gamma$ and let $\gamma\in\Gamma$ satisfy
$\gamma\geq\lambda$ for all $\lambda\in\Gamma$. If $\gamma^r =
\lambda_1\cdots\lambda_r$ for some $r\geq 1$ and
$\lambda_i\in\Gamma$ then $\gamma = \lambda_i$ for $1\leq i \leq
r$.
\end{lem}

\Prf Suppose $\gamma^r = \lambda_1\cdots \lambda_r$ with
$\gamma,\lambda_i\in\Gamma$ for $1\leq i \leq r$ and put
$d(\gamma)=(n_1,\ldots,n_k)$ and
$d(\lambda_i)=(n_1^{(i)},\ldots,n_k^{(i)})$ . Then $r n_j =
\sum_{i=1}^r n_j^{(i)}$ for $1\leq j \leq k$ and since $\gamma
\geq \lambda_i$ we have $n_1 \geq n_1^{(i)} $ for all $i$. This
gives (with $j=1$) $n_1=n_1^{(i)}$ for all $i$. As
$\gamma\geq\lambda_i$, this forces $n_2\geq n_2^{(i)}$ for all $i$
and again (with $j=2$) $n_2 = n_2^{(i)}$ for all $i$. Hence we may
proceed inductively to obtain $n_j=n_j^{(i)}$ for all $1\leq i
\leq r$ and $1\leq j \leq k$. Thus we have $\gamma^r =
\lambda_1\cdots \lambda_r$ with $d(\gamma)=d(\lambda_i)$ for all
$i$, and the result follows from the factorization property. \bx

\begin{thm}\label{radical}
$\fL_\Lambda$ is semisimple if and only if every edge in $\Lambda$
lies on a cycle. If $\Lambda$ has finitely many vertices,
$|\Lambda^0|=n<\infty$, then $\rad(\fL_\Lambda)$ is nilpotent of
degree at most $n$ and is equal to the $\wot$-closed two-sided
ideal generated by $\{L_\lambda: \lambda\in\nc(\Lambda)\}$.
\end{thm}

\Prf Suppose first that every edge in $\Lambda$ lies on a cycle.
Let $A\in\fL_\Lambda$ be non-zero. Then $A\xi_\lambda = R_\lambda
A \xi_{r(\lambda)}$ for all $\lambda$ and so there is some
$v\in\Lambda^0$ such that $A\xi_v = \sum_{s(\lambda)=v} a_\lambda
\xi_\lambda \neq 0$. Let $\Gamma$ be the set of
$\lambda\in\Lambda$ such that $a_\lambda \neq 0$ and
$\delta(\lambda)$ is minimal with this property. Let $\gamma$ be a
maximal element in $\Gamma$ with respect to the lexicographic
ordering discussed above.  By assumption there is a path $\mu$
such that $\mu\gamma$ is a cycle and so the paths $(\mu\gamma)^k$,
$k\geq 1$, are also cycles. Further, $\mu\gamma$ is maximal in the
set $\mu\Gamma$. Thus, by the minimality of $\delta(\gamma)$ and
an application of the previous lemma to $\mu\gamma\in\mu\Gamma$
for each $k\geq 1$, from a consideration of Fourier expansions we
have
\[
(L_\mu A)^k \xi_v = a_\gamma^k \,\,\xi_{(\mu\gamma)^k} +
\sum_{\lambda'\neq (\mu\gamma)^k} b_{\lambda'} \xi_{\lambda'}.
\]
Hence for $k\geq 1$ this yields
\[
||(L_\mu A)^k||^{1 / k} \geq \big| \big< (L_\mu A)^k\xi_v,
\xi_{(\mu\gamma)^k}\big>\big|^{1 / k} = |a_{\gamma}| > 0.
\]
Thus $L_\mu A$ has positive spectral radius and is not
quasinilpotent. Since $0\neq A\in\fL_\Lambda$ was arbitrary it
follows that $\rad(\fL_\Lambda) = \{0\}$ and $\fL_\Lambda$ is
semisimple.

Conversely, suppose that $\lambda\in\nc(\Lambda)\neq \emptyset$.
Since $\lambda$ does not lie on a cycle there are no paths
$\mu_1,\mu_2\in\Lambda$ such that both $\mu_1$ and $\mu_2$ contain
$\lambda$ as an edge and $\mu_1\mu_2$ belongs to $\Lambda$. Hence,
a consideration of Fourier expansions shows that
$(AL_\lambda)^2=0$ for all $A\in\fL_\Lambda$. Thus $L_\lambda$
belongs to $\rad(\fL_\Lambda)$ and $\fL_\Lambda$ has non-zero
radical.

It remains to verify the structure of $\rad(\fL_\Lambda)$ in the
finite vertex case. Let $\J$ be the $\wot$-closed two-sided ideal
in $\fL_\Lambda$ generated by $\{L_\lambda : \lambda\in
\nc(\Lambda)\}$. Suppose first that $A$ belongs to $\rad
(\fL_\Lambda)$ with expansion $A\sim \sum_{\lambda} a_\lambda
L_\lambda$. We claim that a coefficient $a_\lambda$ is non-zero
only if $\lambda$ includes an edge $\lambda'\in\nc(\Lambda)$.
Since the Cesaro sums for $A$ would then belong to $\J$, and they
converge in the strong operator topology to $A$, this would show
that $A$ belongs to $\J$. Suppose by way of contradiction that
there is a path $\lambda$ with $a_\lambda\neq 0$ which includes no
edges from $\nc(\Lambda)$ and, as above, assume $\lambda$ is
maximal in the lexicographic ordering  amongst the paths of
minimal length with this property. Then $\lambda$ belongs to a
transitive component of $\Lambda$. So we may choose $\mu\in
\Lambda$ such that $\mu\lambda$ is a cycle in $\Lambda$ and hence
$(\mu\lambda)^k$ belongs to $\Lambda$ for $k\geq 1$. Then by
assumption $(\mu\lambda)^k$ is a path of minimal length and
satisfies the maximality condition amongst the paths in the
expansion of $(L_\mu A)^k$ that have non-zero coefficients. Hence
an application of the lemma shows  the coefficient of
$L_{(\mu\lambda)^k}$ in this expansion is $(a_\lambda)^k$. Thus we
may argue as above to obtain that $L_\mu A$ has positive spectral
radius and hence is not quasinilpotent. This contradiction
verifies the claim and shows that $\J$ contains
$\rad(\fL_\Lambda)$. Notice that this inclusion does not rely on
finitely many vertices.

For the converse inclusion note that, in the case that
$|\Lambda^0| = n < \infty$, if $\mu_1,\ldots,\mu_n$ are paths in
$\Lambda$, each of which has a factorization that includes at
least one edge from $\nc(\Lambda)$, then the product $\mu_n \cdots
\mu_1$ cannot belong to $\Lambda$. It follows that $\J^n = \{0\}$
in this case because any operator $X$ of the form $X = A_1 \cdots
A_n$ with each $A_i\in\J$ can have no non-zero Fourier
coefficients. Thus, $\J$ is contained in $\rad(\fL_\Lambda)$ and
the result follows.
 \bx

\vspace{0.1in}

{\noindent}{\it Acknowledgements.} The first author was partially
supported by an NSERC grant and an EPSRC visiting fellowship and
the second author was partially supported by an EPSRC grant.

\end{document}